\documentclass{article}
\usepackage{amsmath,amsthm,amssymb,eucal,ifthen, amscd,graphicx,subfigure}
\usepackage{hyperref}

 \numberwithin{equation}{section}

 \newtheorem{theorem}{Theorem}[section]
 \newtheorem{definition}{Definition}[section]

 \newtheorem{proposition}[definition]{Proposition}
 \newtheorem{lemma}[definition]{Lemma}
 \newtheorem{corollary}[definition]{Corollary}
 
 \newtheorem{remark}[definition]{Remark}

 \def\eref#1{(\ref{#1})}
 
 \def\mult#1#2{\textup{mult}_{#1}{(#2)}}

\begin{document}

\title{Vibration Spectra of the $m$-Tree Fractal}

\author{Daniel J. Ford and Benjamin Steinhurst\footnote{Both authors supported in part by the NSF grant DMS-0505622.}}

\maketitle

\begin{center}{\small
Contact:\\
steinhurst@math.uconn.edu\\
Department of Mathematics\\
University of Connecticut\\
Storrs~CT~06269~USA}\end{center}


\begin{abstract}
We introduce a family of post-critically finite fractal trees indexed by the number of branches they possess. Then we produce a Laplacian operator on graph approximations to these fractals and use spectral decimation to describe the spectrum of the Laplacian on these trees. Lastly we consider the behavior of the spectrum as the number of branches increases.

MCS: 28A80,  34B45, 15A18, 60J45, 94C99, 31C25

\end{abstract}

\section{Introduction}\label{sec:intro}
There have been many studies using spectral decimation to calculate the spectrum of a Laplacian operator on finitely ramified fractals \cite{Bajorin2008a,Bajorin2008b,Shima1996,Zhou2008a,Zhou2008b}. These authors have studied many fractals in both specific cases, in families, or in general. We shall consider a Laplacian for a family of fractals, the $m$-Branch Trees and in one general calculation find the spectrum of Laplacians on the entire family. This work is done via the spectral decimation process by which the spectra of Laplacians on graph approximations are used to calculate the spectrum of a limiting Laplacian on the fractal or on an infinite graph. Our notation is slightly complicated by the fact that we do the spectral decimation calculations for an arbitrary $m$, so effectively doing every fractal in this family at one time. Trees have been a well studied topic in fractal literature, for example Vicsek sets in \cite{Zhou2008a, Zhou2008b} and Dendrites in \cite{Kigami1995}.

In \cite{Kigami2001}, Kigami set out the frame work by which post-critically finite fractals can be thought of as abstract spaces independent of any ambient space. This is accomplished by labeling a point by an ``address'' determined by the cell structure of the fractal. The space of addresses is the fractal. Because of this point of view it is natural to consider graph approximations to the fractal as truncations of these addresses. Kigami also proved that the graph Laplacian operators on the approximating graphs converge to a Laplacian on the full fractal. Strichartz, \cite{Strichartz2006}, builds upon Kigami's point of view and considers differential equations on post-critically finite fractals using the types of operators that Kigami constructed. Strichartz included a discussion of the spectral decimation method for calculating the spectrum of these Laplacians from the spectra of the graph Laplacians. 

The  $m$-Branch Trees, $m \geqslant 3$, are post-critically finite fractals whose approximating graphs are defined for $m \geqslant 3$ and are constructed with $m$-simplices. We run through the calculations of the spectrum for the $3$-branch tree in Section\ \ref{sec:three-branch} before doing the general case for the sake of having concrete matrices to work with. Then in Section\ \ref{sec:m-branch} we do the general case and finally in Section\ \ref{sec:infinity} we observe what the behavior of the spectra of the Laplacians on the $m$-branch trees is as the number of branches increase. But first we describe the trees in Section\ \ref{m-tree} and the spectral decimation method in Section\ \ref{sec:spec-dec}.

\section*{Acknowledgments} The authors thank Alexander Teplyaev for his guidance and support, Kevin Romeo for his help during the early stages of this project, and  Denglin Zhou for his insightful comments.

\section{\emph{m}-Branch Tree Fractals}\label{m-tree}
The $m$-branch tree fractal, $F_m$, is a post-critically finite fractal with $m$ defining contraction mappings. The zero-level graph approximation, $V_{m,0}$, consists of a complete graph of $m$ vertices, when $m=3$, $V_{3,0}$ is a triangle, when $m=4$, $V_{4,0}$ is a tetrahedron. The iterated function system that generates the fractal scales, duplicates, and translates the simplex to $m$ simplices sharing a common point at the epicenter of the previous simplex and with each vertex from $V_{m,0}$ as a vertex of one of the new simplices, this is the graph $V_{m,1}$. This process is iterated and the countable set of vertices is completed in the effective resistance metric to form a tree with $m$ branches. With the graphs $V_{m,0}$ and $V_{m,1}$ described we refer to the appendix of \cite{Kigami1993} for the detailed construction of the post-critically finite fractal and proofs of its properties. 

Below we construct a Laplacian operator, $\Delta_m$, on $F_m$ as a limit of Laplacians $\Delta_{m,n}$ on $V_{m,n}$ from which we will calculate the spectrum of $\Delta_m$ by spectral decimation. All of the calculations that are done in the spectral decimation process do not depend on any particular embedding into a Euclidean space and in fact when we calculate the Hausdorff dimension of these space in Section\ \ref{sec:m-branch}\ we will use the intrinsic effective resistance metric. This metric capitalizes on the analogy between graphs and electrical networks to calculate the effective resistance between two points as if the edges of the graph were resistors and the vertices nodes in an electrical network. This analogy has been used by many authors but a useful reference for the mechanics of these calculations and a more extensive bibliography is \cite{Strichartz2006}. It is a straight forward calculation to show that in the effective resistance metric the contraction mappings forming the trees have contraction factor one half for any $m$.

There has recently been developed a connection between fractals and certain groups, most specifically cellular autonoma and iterated monodromy groups. Nekrashevych has written a survey \cite{Nek2005} exploring self-similar groups such as iterated monodromy groups involving Schreier graphs as pictoral representations of the groups. In \cite{BartholdiGrigorchuk2000} and \cite{BartholdiGrigorchukNekrashevych2003} there is consideration of a particular self-similar group, which is a subgroup of the automorphism group of a rooted $p-$tree. This particular self-similar group is sometimes called the Gupta-Fabrykowski group. Bartholdi, Grigorchuk, and Nekrashevych explore the limits of the Schreier graphs as fractals combining algebraic and analytic perspectives. The Schreier graphs for the Gupta-Fabrykowski groups are the same as the $V_{m,n}$ graphs for the $m$-Branch Tree when $m$ is a prime. 

\begin{figure}[htbp]
	\centering
	\begin{tabular}{ccc}
		\includegraphics[height=100pt]{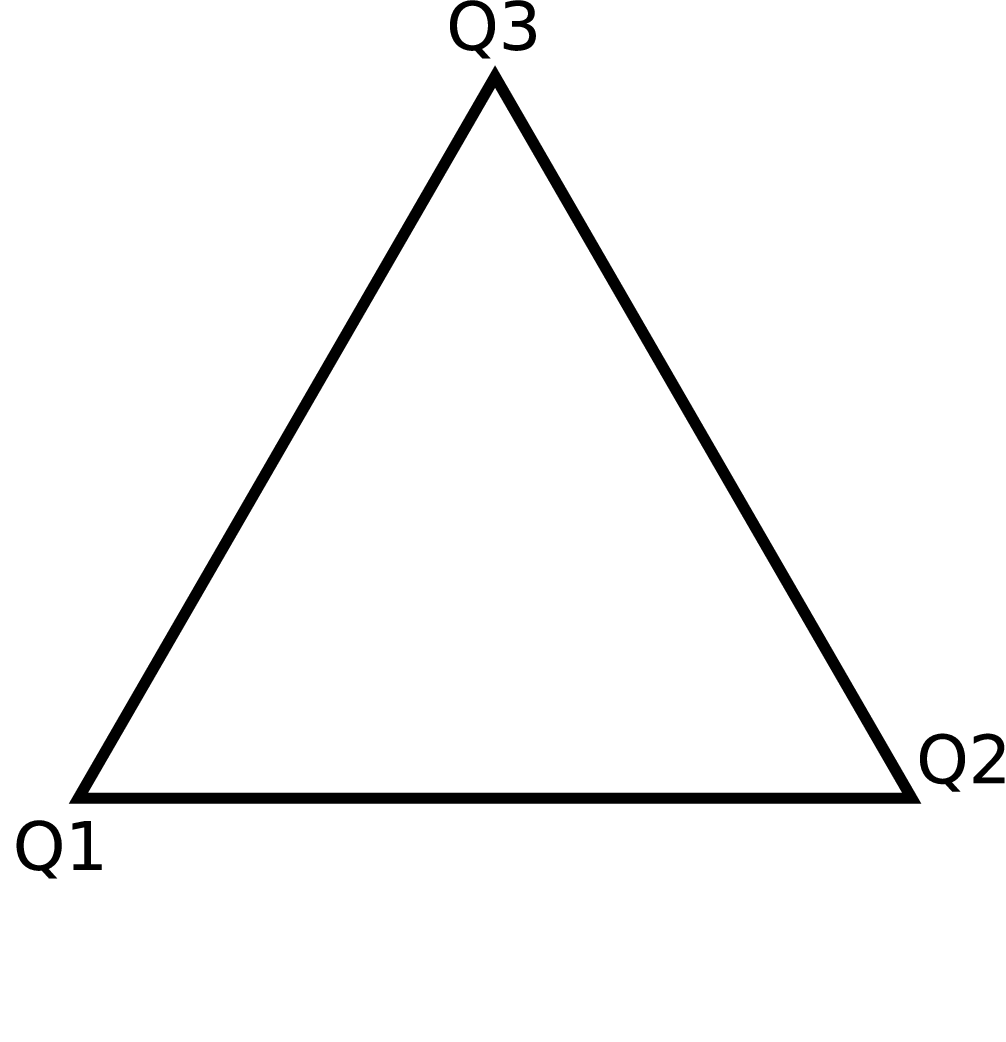} &
		\includegraphics[height=100pt]{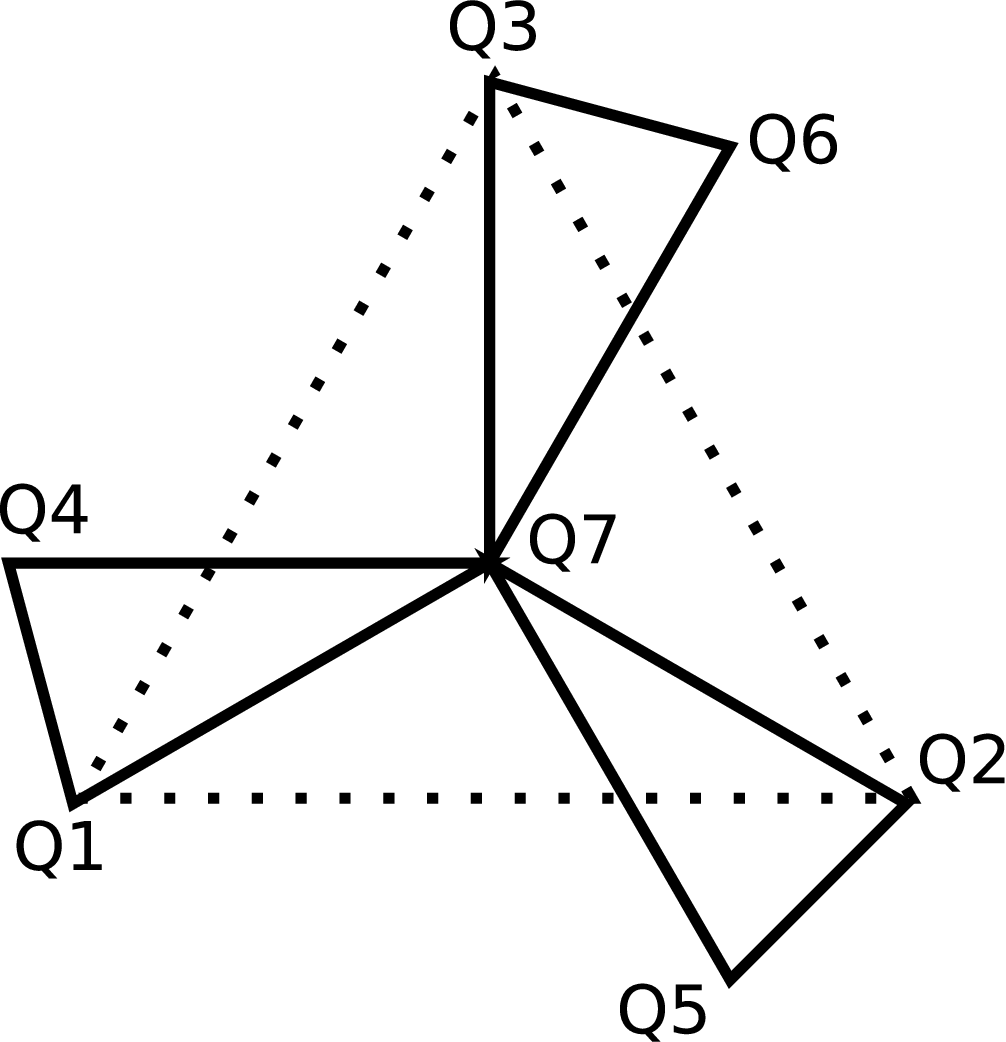} &
		\includegraphics[height=100pt]{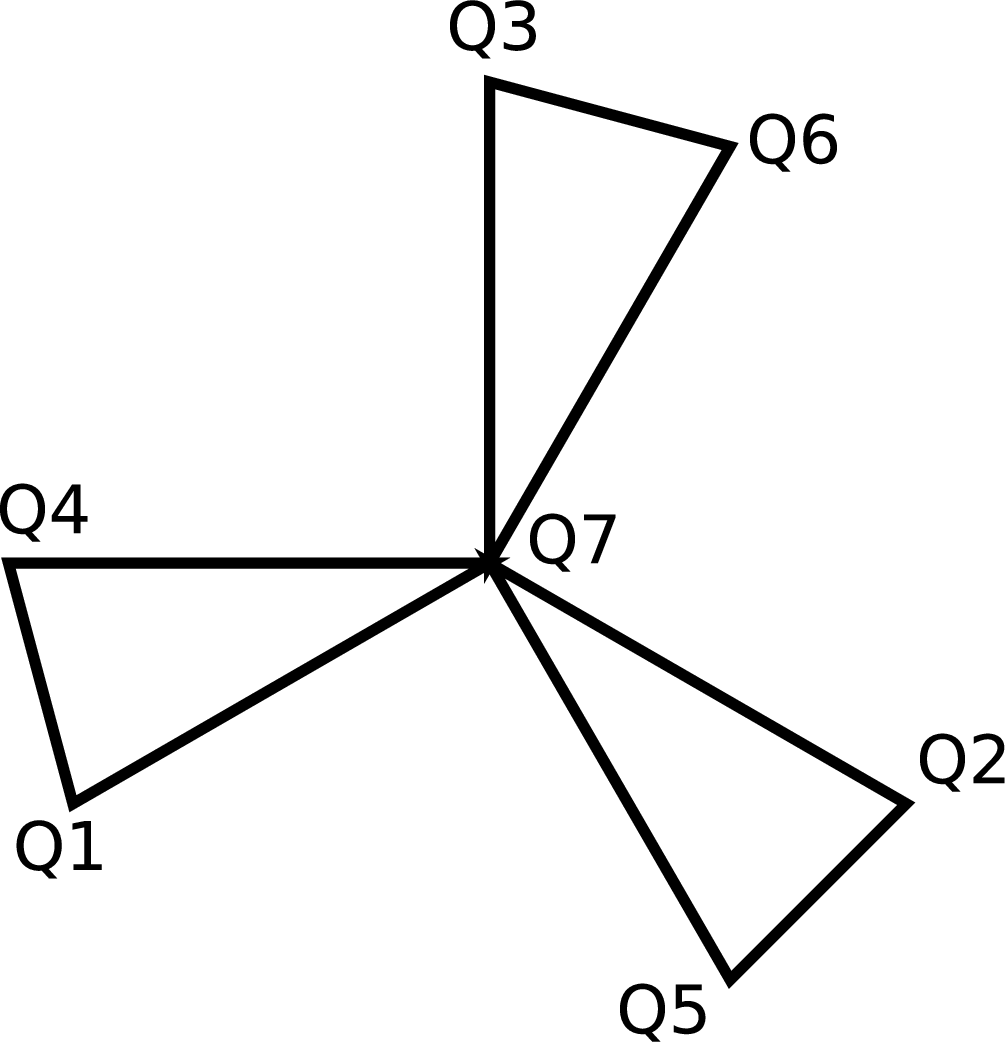} \\
		$V_0$ & $\Rightarrow$ & $V_1$
	\end{tabular}
	\caption{The progression of the three-branch tree from its $V_{3,0}$ to $V_{3,1}$ network}
	\label{fig:F3}
\end{figure}

\begin{figure}[htbp]
	\centering
		\includegraphics[scale=1]{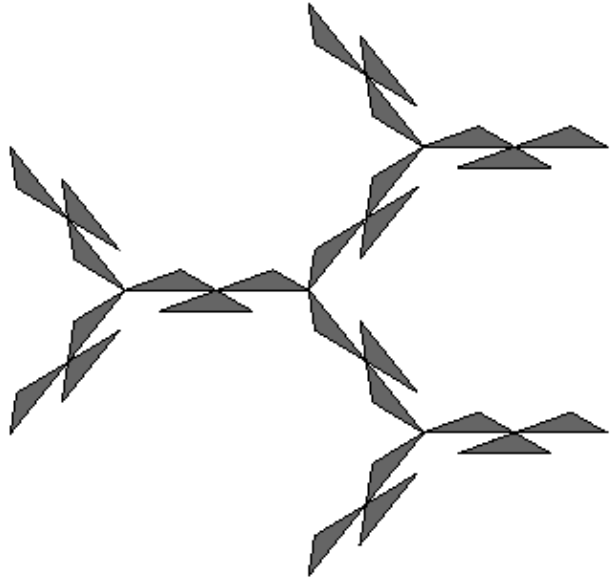}
	\caption{The graph $V_{3,2}$, note that the branches only connect to each other at the center point.}
	\label{fig:F3V2}
\end{figure}

Figures \ref{fig:F3}\ and \ref{fig:F3V2} gives a visual guide to what the iterated function system does in the $m=3$ case. Note that all triangles are equilateral in the intrinsic effective resistance geometry but when embedded in $\mathbb{R}^{2}$ they would overlap if they were drawn as equilateral in the geometry of $\mathbb{R}^{2}$. Since the theoretical machinery that we use in this paper does not depend on a particular embedding into Euclidean space we won't be using the geometry of any $\mathbb{R}^{n}$ so the issue of overlap is not relevant.

The Laplacian operator we use is analogous to the second difference operator on a continuous space because it takes the average value of the function at neighboring points and subtracts that from the value of the function, $\Delta_{m,n} f(x) = f(x)- \sum_{x\sim  y} \frac{f(y)}{\deg(x)}.$ Where the sum is over neighboring points in the graph and $\deg(x)$ is the number of neighbors that $x$ has. Because the matrix representation of this operator on a graph has ones on the diagonal this is known as a probabilistic Laplacian. We construct a matrix representation for the Laplacian $\Delta_{m,n}$ associated to the graph $V_{m,n}$  as in \cite{Bajorin2008a, Bajorin2008b, Strichartz2006, Zhou2008a, Zhou2008b}. Let $M_{m,n}$ denote this matrix which is the  level-$n$ graph approximation to the Laplacian on $F_m$. The on-diagonal entries of $M_{m,n}$ are 1 as already noted. The off-diagonal entries are:
\begin{itemize}
	\item 0 if $x_i$ and $x_j$ are not connected;
	\item $\frac{-1}{deg(x_i)}$ if $x_i$ and $x_j$ \emph{are} connected.
\end{itemize}
Kigami proved that the sequence of operators $\rho^{n}M_n$ converge in an appropriate sense to an operator $\widehat{\Delta}_m$ on $F_m$ in \cite{Kigami2001,Kigami2003} where $\rho$ is the ``energy renormalization constant.'' If we view the distance between adjacent vertices always being constant as the level of approximation increases then we get an infinite graph instead of a compact fractal where $\rho=1$, if the graphs are scaled at each step so that the limiting fractal is bounded then $\rho >1$. We will calculate the spectrum of the bounded operator $\lim_{n\rightarrow \infty} M_{m,n} = \Delta_m$ and will not be concerned with $\rho$.  Since there are $deg(x_i)$ points connected to $x_i$ and the diagonal entry is 1, then we have the property that the row sums of $M_{m,n}$ are always equal to zero which is a very useful aid in the calculations.

\section{Spectral Decimation}\label{sec:spec-dec}
Spectral decimation is a means of extending the eigenvalues for the $n$-level Laplacian to those of the $n+1$-level Laplacian. Used inductively, this means only the zero- and first-level Laplacians must be explicitly calculated. The process, which dates back to \cite{FukushimaShima1992, KigamiLap1993, Shima1996}, is presented in terms of the calculations shown in  \cite{Bajorin2008b} and discussed in more abstract detail in \cite{Bajorin2008a, Shima1996, Strichartz2006}. We first write the level-one Laplacian as a block matrix. For the time being we drop the $m$ from the notation for simplicity's sake.
\begin{equation}\label{abcd}
	M = M_1 = \left(\begin{array}{cc}A & B\\C & D\end{array}\right),
\end{equation}
where $A = I_m$ since the boundary points are not neighbors in $V_1$. The Schur Complement of the matrix $M-Iz$ is the matrix-valued function: 
\begin{equation}\label{schur}
	S(z) = (A-z) - B(D-z)^{-1}C.
\end{equation}
When $v = (v_0,v_1')^{T}$ is an eigenvector of $M$ with corresponding eigenvalue $z$, the eigenvalue equation can be written as
\begin{equation}\label{eigen}
	\left(\begin{array}{cc}A&B\\C&D\end{array}\right) \left(\begin{array}{cc}v_0\\v_1^{\prime}\end{array}\right)
		= z \left(\begin{array}{cc}v_0\\v_1^{\prime}\end{array}\right).
\end{equation}
This system of linear equations yields two equations: $v_1^{\prime}=-(D-z)^{-1}Cv_0$, provided that $z \notin \sigma(D)$, and $(A-z)v_0 + Bv_1^{\prime}=0$. It is worth noting that the map $-(D-z)^{-1}C$ takes the eigenvector $v_0$ on the boundary points and determines the values required at the interior vertices so that $\left(\begin{array}{cc}v_0\\v_1^{\prime}\end{array}\right)$ is an eigenvector on $V_1$. This is what actually does the eigenvector extension mentioned below. When the two equations are combined with the Schur complement they imply that
\begin{equation}\label{Sziszero}
	S(z)v_0 = 0.
\end{equation}
If $v_0$ is \emph{also} an eigenvector of $M_0$ with corresponding eigenvalue $z_0$, then
\begin{equation}\label{eigenM0}
	(M_0 - z_0)v_0 = 0.
\end{equation}
By combining \eref{Sziszero} and \eref{eigenM0}, and allowing $z_0 = R(z)$ for a rational function $R(z)$, we have,
\begin{equation}\label{final}
	S(z) = \phi(z) \left(M_0 - R(z) \right)
\end{equation}
where $\phi(z)$ and $R(z)$ are scalar-valued rational functions whose existence is proved in \cite{Bajorin2008a, KigamiLap1993, Shima1996, Strichartz2006}. It can be observed from the proof of their existence that $\phi(z)$ is dependent on $S_{1,2}$ and $R(z)$ is dependent on both $S_{1,1}$ and $\phi(z)$. The role of $R(z)$ is as an ``eigenvalue projector'' taking eigenvalues of one level and projecting them down to a lower level's eigenvalues. To calculate these functions we only need two entries from $S(z)$. The formulae are: 
\begin{equation}
	\phi(z) = -(m-1)S_{1,2} \hspace{1cm} and \hspace{1cm} R(z) = 1-\frac{S_{1,1}}{\phi(z)}.
\end{equation}

These equations are \emph{a priori} rational functions on $\mathbb{C}$ so the locations of possible zeroes and poles require special attention. If they are off  the real axis they won't be an issue since $\Delta_m$ are bounded symmetric operators, only zeros and poles on the real axis are of concern. The location of these possible zeroes and poles are what we call exceptional values. The exceptional values are $E(M_0,M) = \sigma(D) \cup \left\{z : \phi(z) = 0\right\}$. These values are exceptional because $(D-z)$ is not invertible when $z \in \sigma(D)$, which causes \eref{Sziszero} to not be defined, and if $\phi(z) = 0$ then either \eref{Sziszero} and/or \eref{eigenM0} fail to be defined, i.e. the poles of $R$ and $\phi$ and zeroes of $\phi$. 

It is the rational functions $R(z)$ and $\phi(z)$ that are the main tools used throughout this paper. Proposition \ref{prop:schur} will produce a general form for the values of $S_{1,1}$ and $S_{1,2}$ to obtain the functions $\phi(z)$ (Corollary \ref{c:phi}) and $R(z)$ (Corollary \ref{c:R}) in the general $m$-branch case. The function $R(z)$ projects level-$n+1$ eigenvalues to level-$n$ eigenvalues, so the primary use of $R(z)$ is to take inverse images under it. We will call the depth-$n$ eigenvalue $z$ the offspring of the depth-$n+1$ eigenvalue $R^{-1}(z)$.

The process of extending eigenvalues from one approximation to the next is summarized in the following proposition from \cite{Bajorin2008a}. Above an extension map was mentioned, $-(D-z)^{-1}C$, from this extension map the eigenvectors are extended providing that $z$ is not an exceptional value but is an eigenvalue at the next higher level. To find those eigenvalues at the next level we will use this proposition. Denote by $\textup{mult}_D(z)$ the multiplicity of $z$ as an eigenvalue of $D$. Similarly, $\textup{mult}_n(z)$ is the multiplicity of $z$ as an eigenvalue of $M_n$. Again the subscript $m$ has been removed to simplify notation, so everything here is for a  given value of $m$. 

\begin{remark}\label{remark:dim}
Here and throughout, $\textup{dim}_n$ is the dimension of the function space on $V_n$ for a given $m$. Since $V_n$ is a finite collection of points, the dimension of the space of functions on $V_n$ is just the number of points in $V_n$.
\end{remark}

\begin{proposition}\label{p-mult}\cite{Bajorin2008a}
\begin{enumerate}
        \item\label{p-mult1}
If $z\notin E(M_0, M)$, then
\begin{equation}
	\mult n{z}=\mult{n-1}{ R(z)},
\end{equation}
and every corresponding eigenfunction at depth $n$ is an extension of an eigenfunction at depth $n-1$.

        \item\label{p-mult2}
If $z \notin \sigma(D)$, $\phi(z) = 0$ and $R(z)$ has a removable singularity at $z$, then
\begin{equation}
	\mult n{z}=\textup{dim}_{n-1},
\end{equation}
and every corresponding eigenfunction at depth $n$ is localized.

        \item\label{p-mult3}
If $z \in \sigma(D)$, both $\phi(z)$ and $\phi(z)R(z)$ have poles at $z$, $R(z)$ has a removable singularity at $z$,  and $\frac d{dz}R(z)\neq 0$, then
\begin{equation}
	\mult n{z}=m^{n-1}\mult D{z}-\textup{dim}_{n-1}+\, \mult{n-1}{ R(z)},
\end{equation}
and every corresponding eigenfunction at depth $n$ vanishes on $V_{n-1}$.

        \item\label{p-mult4}
If $z \in \sigma(D)$, but $\phi(z)$ and $\phi(z)R(z)$ do not have poles at $z$, and $\phi(z)\neq 0$, then
\begin{equation}
	\mult n{z}=m^{n-1}\mult D{z}+\, \mult{n-1}{R(z)}.
\end{equation}
In this case $m^{n-1}\mult D{z}$ linearly independent eigenfunctions are localized, and $\mult{n-1}{ R(z)}$ more linearly independent eigenfunctions are extensions of corresponding eigenfunction at depth $n-1$.

        \item\label{p-mult5}
If $z \in \sigma(D)$, but $\phi(z)$ and $\phi(z)R(z)$ do not have poles at $z$, and $\phi(z)= 0$, then
\begin{equation}
	\mult n{z}=m^{n-1}\mult D{z}+\, \mult{n-1}{ R(z)}+\, \textup{dim}_{n-1}
\end{equation}
provided $R(z)$ has a removable singularity at $z$. In this case there are $m^{n-1}\mult D{z}+\textup{dim}_{n-1}$ localized and
$\mult{n-1}{ R(z)}$ non-localized corresponding eigenfunctions at depth $n$.

        \item\label{p-mult6}
If $z \in \sigma(D)$, both $\phi (z)$ and $\phi(z)R(z)$ have poles at $z$, $R(z)$ has a removable singularity at $z$, and $\frac d{dz}R(z)= 0$, then
\begin{equation}
	\mult n{z}=\mult{n-1}{ R(z)},
\end{equation}
provided there are no corresponding eigenfunctions at depth $n$ that vanish on $V_{n-1}$. In general we have
\begin{equation}
	\mult n{z}=m^{n-1}\mult D{z}-\textup{dim}_{n-1}+\, 2\mult{n-1}{ R(z)}.
\end{equation}

        \item\label{p-mult7}
If $z \notin \sigma(D)$, $\phi(z)=0$ and $R(z)$ has a pole $z$, then $\mult n{z}=0$ and $z$ is not an eigenvalue.

        \item\label{p-mult8}
If $z \in \sigma(D)$, but $\phi(z)$ and $\phi(z)R(z)$ do not have poles at $z$, $\phi(z)= 0$, and $R(z)$ has a pole $z$, then
\begin{equation}
	\mult n{z}=m^{n-1}\mult D{z}
\end{equation}
and every corresponding eigenfunction at depth $n$ vanishes on $V_{n-1}$.
\end{enumerate}
\end{proposition}

	\def\pro#1{Proposition~\ref{p-mult}(\ref{p-mult#1})}
	\def\Pro{Proposition~\ref{p-mult}}

The proof of this proposition can be found in \cite{Bajorin2008a} and in a large part depends on the Schur Complement formula and the eigenvector equations from previously in the section. In this source it is shown how the eigenspace projection matrices are constructed, these matrices are related to the map $-(D-z)^{-1}C$ and project eigenvectors on $V_n$ onto $V_{n-1}$. This proposition is used constantly as we recursively extend eigenvalues from the $V_1$ approximation, which are calculated explicitly, to the full fractal $F_m$. A necessary condition for $\sigma (M_{m,n})$ capturing all of $\sigma(\Delta_m)$ is that the spectral dimension is less than two \cite{Shima1996} however it is known that for p.c.f. fractals this is always the case.

\section{The Three-Branch Tree}\label{sec:three-branch}
We analyze the spectrum of the Laplacian, $\Delta_3$, on the three-branch tree as a special case to demonstrate the calculations involved with spectral decimation as found in \cite{Bajorin2008a, Bajorin2008b} while the approximating graphs and matrices are still relatively small. For this section assume that $m=3$ throughout with will drop the $m$ from the notation. Spectral decimation uses the spectra of $M_0$, $M_1$, $D$, and the Schur Complement of $M_1 - Iz$ yielding the spectrum of $M_n$ as the output. We will then take the limit as $n$ grows to infinity to find $\sigma(\Delta_3)$. The process will also be used in the next section as we consider the general case of the $m$-branch tree.

First we consider $V_0$, which is simply a triangle (Figure \ref{fig:F3}). The Laplacian matrix for $V_0$ is
$$M_0 = \left(\begin{array}{rrr}
	1 & -\frac{1}{2} & -\frac{1}{2}\\
	-\frac{1}{2} & 1 & -\frac{1}{2}\\
	-\frac{1}{2} & -\frac{1}{2} & 1
\end{array}\right).$$
Since $V_0$ is a complete graph there are no zeros in this matrix.

The eigenvalues of $M_0$ written with multiplicities are: $$\sigma(M_0) = \left\{0, \frac{3}{2}, \frac{3}{2}\right\}.$$

Next we consider the Laplacian on $V_1$, which is a set of three triangles with a common vertex (Figure \ref{fig:F3}). There are seven vertices in $V_1$, and the depth-one Laplacian matrix is:
$$M_1 = \left( \begin{array}{rrrrrrr}
	1 & 0 & 0 & -\frac{1}{2} & 0 & 0 & -\frac{1}{2} \\
	0 & 1 & 0 & 0 & -\frac{1}{2} & 0 & -\frac{1}{2} \\
	0 & 0 & 1 & 0 & 0 & -\frac{1}{2} & -\frac{1}{2} \\
	-\frac{1}{2} & 0 & 0 & 1 & 0 & 0 & -\frac{1}{2} \\
	0 & -\frac{1}{2} & 0 & 0 & 1 & 0 & -\frac{1}{2} \\
	0 & 0 & -\frac{1}{2} & 0 & 0 & 1 & -\frac{1}{2} \\
	-\frac{1}{6} & -\frac{1}{6} & -\frac{1}{6} & -\frac{1}{6} & -\frac{1}{6} & -\frac{1}{6} & 1
\end{array}\right).$$
Recall from \eref{abcd} that $A$ will be the identity matrix $I_3$. 

The Schur Complement is $S(z) = (A-z) - B(D-z)^{-1}C$, as stated earlier. With these small matrices $S(z)$ can easily be computed:
$$S(z) = \left(\begin{array}{ccc}
 -\frac{(2z-3)(6z^2-6z+1)}{6(z-1)(2z-1)} & \frac{2z-3}{12(z-1)(2z-1)} & \frac{2z-3}{12(z-1)(2z-1)}\\
 \frac{2z-3}{12(z-1)(2z-1)} & -\frac{(2z-3)(6z^2-6z+1)}{6(z-1)(2z-1)} & \frac{2z-3}{12(z-1)(2z-1)}\\
 \frac{2z-3}{12(z-1)(2z-1)} & \frac{2z-3}{12(z-1)(2z-1)} & -\frac{(2z-3)(6z^2-6z+1)}{6(z-1)(2z-1)}
\end{array}\right).$$
From this, we have
	$$\phi(z) = -(m-1)S_{1,2} = \frac{3-2z}{6(-1+z)(-1+2z)}$$ and 
	$$R(z) = 1 - \frac{S_{1,1}}{\phi(z)} = 6z-6z^2,$$
as defined in \cite{Bajorin2008a, Bajorin2008b}. These two functions allow us to use Proposition \ref{p-mult} to determine the multiplicities of eigenvalues for $M_n$, without writing the matrix $M_n$ down. These recursive formulae utilize the spectrum of $M_1$ to recursively list all eigenvalues of $M_n$.

The eigenfunction extension map that fills in the interior values when eigenvectors are extended is:
$$(D-z)^{-1}C = \left(\begin{array}{cccc}
	\frac{-2+3z}{3(-1+z)(-1+2z)} & \frac{-1}{6(-1+z)(-1+2z)} & \frac{-1}{6(-1+z)(-1+2z)} \\
	\frac{-1}{6(-1+z)(-1+2z)} & \frac{-2+3z}{3(-1+z)(-1+2z)} & \frac{-1}{6(-1+z)(-1+2z)} \\
	\frac{-1}{6(-1+z)(-1+2z)} & \frac{-1}{6(-1+z)(-1+2z)} & \frac{-2+3z}{3(-1+z)(-1+2z)} \\
	\frac{1}{3(-1+2z)} & \frac{1}{3(-1+2z)} & \frac{1}{3(-1+2z)}
\end{array}\right).$$

\begin{figure}[htbp]
	\centering
		\includegraphics[height=100pt]{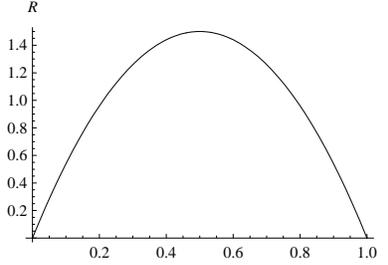}
	\caption{The graph of $R(z)$ for the three-branch tree}
	\label{fig:R3}
\end{figure}

Now we list the eigenvalues of $M_1$, written with multiplicities, to provide a check for the calculations.
	$$\sigma(M_1) = \left\{\frac{3}{2}, \frac{3}{2}, \frac{3}{2}, \frac{3}{2}, \frac{1}{2}, \frac{1}{2}, 0\right\}.$$
With corresponding eigenvectors of \{-1, -1, -1, 0, 0, 0, 1\}, \{0, 0, -1, 0, 0, 1, 0\}, \{0, -1, 0, 0, 1, 0, 0\}, \{-1, 0, 0, 1, 0, 0, 0\}, \{-1, 0, 1, -1, 0, 1, 0\}, \{-1, 1, 0, -1, 1, 0, 0\}, and \{1, 1, 1, 1, 1, 1, 1\} respectively.

To find the exceptional set we need the eigenvalues of $D$, which are written with multiplicities:
	$$\sigma(D) = \left\{\frac{3}{2}, 1, 1, \frac{1}{2}\right\},$$
with corresponding eigenvectors of \{-1, 1, -1, 1\}, \{-1, 0, 1, 0\}, \{-1, 1, 0, 0\}, and \{1, 1, 1, 1\} respectively.

The equation $\phi(z) = 0$ has only one solution, $\left\{\frac{3}{2}\right\}$, so the exceptional set is:
	$$E(M_0,M_1) = \sigma(D) \cup \{z : \phi(z) = 0\} = \left\{\frac{3}{2}, 1, \frac{1}{2}\right\}.$$

We can find the multiplicities of these exceptional values by using \Pro. For the value of $z = \frac{3}{2}$, which is in $\sigma(D)$, is not a pole of $\phi(z)$, and $\phi(z) = 0$, we use \pro5 to find the multiplicities:
$$\begin{array}{l}
	\textup{mult}_0(\frac{3}{2}) = 2,\\
	\textup{mult}_1(\frac{3}{2}) = 4,\\
	\textup{mult}_2(\frac{3}{2}) = 10.
\end{array}$$

For the values of $z = \frac{1}{2}$ and $z = 1$, which are poles of $\phi(z)$ and in $\sigma(D)$, we use \pro3:
$$\begin{array}{ll}
	\textup{mult}_0(\frac{1}{2}) = 0, & \textup{mult}_0(1) = 0,\\
	\textup{mult}_1(\frac{1}{2}) = 2, & \textup{mult}_1(1) = 0,\\
	\textup{mult}_2(\frac{1}{2}) = 4, & \textup{mult}_2(1) = 0.
\end{array}$$

For the value of $z = 0$, since $0 \notin E(M_0,M)$, we use \pro1 to determine that $\textup{mult}_n(0) = 1$.

Using the inverse function $R^{-1}(z)$, we can calculate the ancestors of $z = \frac{1}{2}$ to be $R^{-1}(1/2) = \frac{3\pm\sqrt{6}}{6}$. From \pro1, since $\frac{3\pm\sqrt{6}}{6} \notin E(M_0,M_1)$, their multiplicity at depth $n$ will be the multiplicity of $R(\frac{3\pm\sqrt{6}}{6}) = \frac{1}{2}$ at depth $n-1$.

Table \ref{tab:ancestor3} shows the ancestor-offspring structure of the eigenvalues of the three-branch tree. The symbol * indicates an ancestor of $\frac{3\pm\sqrt{6}}{6}$, calculated by the inverse function $R^{-1}(z)$ computed at the ancestor value $z$. By \pro1 the ancestor and the offspring have the same multiplicity. The empty columns represent exceptional values. If they are eigenvalues of the appropriate $M_n$, then their multiplicity is shown in the right hand part of the same row.

	\global\newcount\NN
	\def\newz#1{\global\NN=0 $z\in\sigma(M_{#1})$\vbox to 4ex{}}
	\def\newm#1{\global\NN=1 $\textup{mult}_{#1}(z)$\vbox to 4ex{}}
	\def\m#1#2#3{\multicolumn{#1}{c|}{\text{\hskip-.2em\ifcase\NN
	\mathversion{bold}$#2$\or\bfseries\itshape#3\fi\hskip-.2em}}}
\begin{table}[htbp]
	\def\newz#1{\global\NN=0 \text{\small{\hskip-.3em}$z{\in}\sigma(M_{#1})$\vbox to 4ex{}{\hskip-.5em}}}
	\def\newm#1{\small\global\NN=1 {\hskip-.3em}$\textup{mult}_{#1}(z)${\hskip-.3em}\vbox to 4ex{}}
	\def\m#1#2#3{\multicolumn{#1}{c|}{\text{\small\hskip-.3em\ifcase\NN
	\small$#2$\or\bfseries\itshape\small#3\fi\hskip-.3em}}}%
        \centering\small
\begin{tabular}{|l||*{22}{c}}
	\cline{1-9}
		\newz0 &\m401 &\m4{\frac32}2 \\
		\newm0 &\m401 &\m4{\frac32}2 \\
		[2ex]\cline{1-13}
		\newz1 &\m301 &\m1{1}{} &\m4{\frac12}2 &\m4{\frac32}4 \\
		\newm1 &\m301 &\m1{}{}  &\m4{\frac12}2 &\m4{\frac32}4 \\
		[2ex]\cline{1-4}\cline{6-15}
		\newz2 &\m201 &\m1{1}{} &\m1{}{} &\m4{\frac{3\pm\sqrt{6}}{6}}{2} &\m4{\frac12}{4} &\m2{\frac32}{10}\\
		\newm2 &\m201 &\m1{1}{} &\m1{}{} &\m2{\frac{3+\sqrt{6}}{6}}{2} 
			&\m2{\frac{3-\sqrt{6}}{6}}{2} &\m4{\frac12}{4} &\m2{\frac32}{10}\\
		[2ex]\cline{1-3}\cline{6-16}
		\newz3 &\m101 &\m1{1}{} &\m1{}{} &\m1{}{} &\m1{*}{2} &\m1{*}{2} &\m1{*}{2} &\m1{*}{2}
			&\m4{\frac{3\pm\sqrt{6}}{6}}{4} &\m2{\frac12}{10} &\m1{\frac32}{28}\\
		\newm3 &\m101 &\m1{1}{} &\m1{}{} &\m1{}{} &\m1{*}{2} &\m1{*}{2} &\m1{*}{2} &\m1{*}{2}
			&\m2{\frac{3+\sqrt{6}}{6}}{4} &\m2{\frac{3-\sqrt{6}}{6}}{4} &\m2{\frac12}{10} &\m1{\frac32}{28}\\%
		[2ex]\cline{1-16}
\end{tabular}\vskip2ex
	\caption{Ancestor-offspring structure of the eigenvalues on the three-branch tree}
	\label{tab:ancestor3}
\end{table}

\begin{proposition}\label{summary-3}
The dimension of the function space on $V_n$, $\dim_n$, as was commented in Remark \ref{remark:dim}, is the same as the number of points in $V_n$. For the three-branch tree, this is:
	$$\textup{dim}_n = 1 + 2\cdot3^n.$$
Likewise, $\textup{mult}_n(z)$ is given as follows for $n \geqslant 0$ and $1 \leqslant k \leqslant n$:
$$\begin{array}{l}
	\textup{mult}_n(0) = 1,\\
	\textup{mult}_n(1) = 0,\\
	\textup{mult}_n(\frac{3}{2}) = 1 + 3^n,\\
	\textup{mult}_n(\frac{1}{2}) = 1 + 3^{n-1},\\
	\textup{mult}_n(\frac{3\pm\sqrt{6}}{6}) = 1 + 3^{n-2},\\
	\indent \vdots\\
	\textup{mult}_n(R^{-k}(\frac{3}{2})) = 1 + 3^{n-k}.
\end{array}$$
\end{proposition}
\begin{proof}
For $F_3$, the number of vertices begins with 3 in $V_0$ then increases by $3^{n-1} + 3^n$ for each subsequent level. Therefore, the partial sum of this sequence yields $\textup{dim}_n = 1 + 2\cdot3^n$. In Proposition \ref{prop:dim}, we will prove $\textup{dim}_n$ for the general case.

By applying Proposition \ref{p-mult} to the eigenvalues of $M_1$, and by using $\textup{dim}_n = 1 + 2\cdot3^n$, $\phi(z)$, and $R(z)$, the multiplicities of $\sigma(M_n)$ are calculated inductively.
\end{proof}

	\def\Julia#1#2{\mathcal{J}_{#1}{(#2)}}
\begin{lemma}\label{julia-3}
The spectrum of $\Delta_3$ is $\bigcup_{i \ge 0}\left\{ R^{-i}\left(\frac{3}{2}\right)\right\} \cup \{0\}$. This spectrum is bounded, $\sigma(\Delta_3) \subseteq [0,1) \cup \{\frac{3}{2}\}$, and accumulates to the Julia set of $R(z)$.
\end{lemma}
\begin{proof}
From Figure \ref{fig:R3} we can see that $R(z)$ is a parabola with zeros $\{0,1\}$ and vertex $(\frac{1}{2}, \frac{3}{2})$. By applying inverses of $R(z)$, we have $R^{-1}(z) \in [0,1]$ so long as $z \in [0,\frac{3}{2}]$. Since $R^{-1}(0) = \{0,1\}$ and $\textup{mult}_n(1) = 0$, the only ancestor of 0 will be itself. Since $\bigcup_{i \ge 0}\left\{ R^{-i}\left(\frac{3}{2}\right)\right\} \subseteq (0,1) \cup \{\frac{3}{2}\}$, we have
$$\sigma(\Delta_3) = \bigcup_{i \ge 0}\left\{ R^{-i}(\tfrac{3}{2})\right\} \cup \{0\} \subseteq [0,1) \cup \{\tfrac{3}{2}\}.$$
The statement about the Julia set is a special case of the discussion in \cite{Falconer1990} Chapter 14.
\end{proof}

\section{The \emph{m}-Branch Tree}\label{sec:m-branch}
This section focuses on the use of spectral decimation on the $m$-branch tree and its Laplacian, $\Delta_m$. We follow the same process as in Section \ref{sec:three-branch}, however the inverse of $(D-z)$ is more difficult to compute since it must be done for arbitrary $m$ so we must use a different method for this computation.

\begin{proposition}\label{prop:dim}
For $F_m$, the number of vertices in $V_{m,n}$, which is equal to the dimension of the function space on $V_{m,n}$ (Remark \ref{remark:dim}), is
$$\textup{dim}_{m,n} = 1 + (m-1)m^n,$$
where $m$ and $n$ are integers with $m \geqslant 3$ and $n \geqslant 0$.
\end{proposition}
\begin{proof}
Beginning with $V_{m,0}$ in $F_m$, $dim_{m,0} = m$. At each subsequent level, the number of vertices introduced by the defining iterated function system is $(m-2)m^n + m^{n-1}$. Thus,
	$$\textup{dim}_{m,n} = m + \sum^{n}_{i=1}(m-2)m^i + m^{i-1}.$$
For any level $n \geqslant 0$, the partial sum of the sequence is $\textup{dim}_n = 1 + (m-1)m^n$.
\end{proof}

\begin{proposition}\label{prop:mult0}
For the depth-zero Laplacian matrix, $\textup{mult}_0(0) = 1$ and $\textup{mult}_0(\frac{m}{m-1}) = m-1$.
\end{proposition}
\begin{proof}
For $F_m$, we will have an $m \times m$ depth-zero Laplacian matrix of the form
$$M_{m,0} = \left(\begin{array}{ccccc}
	1 & \frac{-1}{m-1} & \frac{-1}{m-1} & \cdots & \frac{-1}{m-1}\\
	\frac{-1}{m-1} & 1 & \frac{-1}{m-1} & \cdots & \frac{-1}{m-1}\\
	\frac{-1}{m-1} & \frac{-1}{m-1} & 1 & \cdots & \frac{-1}{m-1}\\
	\vdots & & & \ddots & \vdots\\
	\frac{-1}{m-1} & \frac{-1}{m-1} & \frac{-1}{m-1} & \cdots & 1
\end{array}\right).$$

We are looking for the values of $z$ and vectors, $v$, such that $M_{m,0} v = z v$. Since $M_{m,0}$ has a row-sum of zero $v = \{1, \ldots, 1\} \in \mathbb{R}^m$ will be an eigenvector of $M_{m,0}$ with corresponding eigenvalue $z = 0$. So this is our first pair.

Consider $v = \{0, \ldots, 0, 1, -1, 0, \ldots, 0 \} \in \mathbb{R}^m$, where 1 and -1 occur consecutively and all other elements of $v$ are 0, we get a corresponding eigenvalue $z = \frac{m}{m-1}$. There are $m-1$ of these vectors, so $\textup{mult}_0(\frac{m}{m-1}) = m-1$.

These vectors are linearly independent by standard methods, so they form a full set of eigenvectors for $M_{m,0}$.
\end{proof}

For the $m$-branch tree, $M_{m,1}$ is a square matrix with $m^2-m+1$ rows and columns. Take $A$ to be the $m \times m$ identity matrix, then $D$ is a square matrix with $(m-1)^2$ rows and columns. Likewise, $B$ will have dimension $m \times (m-1)^2$ and $C$ will have dimension $(m-1)^2 \times m$. Let
$$\delta = \left(\begin{array}{cccc}
	\frac{-1}{m-1} & 0 & \cdots & 0\\
	0 & \frac{-1}{m-1} & \cdots & 0\\
	\vdots & & \ddots & \vdots\\
	0 & 0 & \cdots & \frac{-1}{m-1}
\end{array}\right)_{m \times m}	= \frac{-1}{m-1}I_m.$$
We will then have:
$$M_{m,1} = \left(\begin{array}{c|c}
	A & B\\
	\hline
	C & D
\end{array}\right) = 
\left(\begin{array}{c|cccc}
	I & \delta & \cdots & \delta & \frac{-1}{m-1}\\
	\hline
	\delta & I & \cdots & \delta & \frac{-1}{m-1}\\
	\vdots &  & \ddots & & \vdots\\
	\delta & \delta & \cdots & I & \frac{-1}{m-1}\\
	\frac{-1}{m(m-1)} & \frac{-1}{m(m-1)} & \cdots & \frac{-1}{m(m-1)} & 1
\end{array}\right).$$
In the last row and column the entries are row and column vectors with the indicated values in every entry. The lower right hand entry is a single entry however. 

\begin{proposition}\label{prop:sigmaD} For the $m$-branch tree, $\sigma(D) = \left\{\frac{m}{m-1}, \frac{2}{m-1}, \frac{1}{m-1}\right\}$ with multiplicities $\textup{mult}_D(\frac{m}{m-1}) = m^2-3m+1, \textup{mult}_D(\frac{2}{m-1}) = m-1,$ and $\textup{mult}_D(\frac{1}{m-1}) = 1$.
\end{proposition}

\begin{proof}
We can observe that $D$ has a row-sum of $\frac{1}{m-1}.$ Therefore, $$v_1 = \{1, \ldots, 1\}$$ will be an eigenvector of $D$ with corresponding eigenvalue $\lambda_1 = \frac{1}{m-1}$.

We state the remaining eigenvectors and their eigenvalues. If $$v_{2a} = \{-1, \ldots, -1, 0, \ldots, 0, 1\},$$ where -1 occurs in the first $m$ entries of the vector, $v_{2a}$ will be an eigenvector of $D$ with eigenvalue $\lambda_2 = \frac{m}{m-1}$.

Consider $$v_{2b} = \{x_1^1, x_2^1, \ldots, x_m^1, x_1^2, \ldots, x_m^{m-2}, 0\},$$ where $\sum_{i=1}^{m-2} x_j^i = 0$, note that the sum is across the upper index. If we choose $x_1^1 = -1$, for some $i >1$, $x_i^1 = 1$, and all other $x_j = 0$, $v_{2b}$ will be an eigenvector of $D$ with eigenvalue $\lambda_2 = \frac{m}{m-1}$. Since the 1 can occur in one of $m$ entries, and the corresponding -1 can occur in one of $m-3$ entries, $v_{2b}$ occurs in $m(m-3) = m^2-3m$ linearly independent variations, and consequently $\lambda$ has a multiplicity of $m^2-3m+1$ (including the vector $v_{2a}$).
	\def\vx{\{x_1^1, x_2^1, \ldots, x_m^1, x_1^2, \ldots, x_m^{m-2}, 0\}}
 
Consider $$v_3 = \{c_1, c_2, \ldots, c_m, c_1, \ldots, c_m, \ldots, 0\},$$ where $\sum_{i=1}^{m}c_i = 0$. If we choose $c_1 = -1$, some $c_i = 1$, and all other $c_j = 0$, $v_3$ will be an eigenvector of $D$ with eigenvalue $\lambda_3 = \frac{2}{m-1}$. Since there are $m-1$ of these vectors, $\lambda$ has a multiplicity of $m-1$.
	\def\vc{\{c_1, c_2, \ldots, c_m, c_1, \ldots, c_m, \ldots, 0\}}
	
It can be checked that $v_1$, $v_{2a}$, $v_{2b}$, and $v_3$ are linearly independent, so these vectors form a spanning set of eigenvectors.
\end{proof}

\begin{proposition}\label{prop:schur}
The two entries of the Schur complement used in calculating $R(z)$ and $\phi(z)$ are:
	$$S_{1,1} = (1-z) - \frac{m-z(m-1)^2}{m(1-z(m-1))(2-z(m-1))}$$ and
	$$S_{1,2} = -\frac{m-z(m-1)}{m(m-1)(1-z(m-1))(2-z(m-1))}$$
for any $m \geqslant 3$.
\end{proposition}

\begin{proof}
To find the Schur Complement, one must invert the increasingly large matrix $D-z$, where $D$ is an $(m-1)^2 \times (m-1)^2$ matrix. This is most conveniently done by forming a spectral resolution of $D$. Since $D$ has 3 distinct eigenvalues, we can express $D$ as:
	$$D = \lambda_1 \cdot P_1 + \lambda_2 \cdot P_2 + \lambda_3 \cdot P_3,$$
where $\lambda_i$ is an eigenvalue of $D$ and $P_i$ is the projection matrix onto the eigenspace, $E_i$, corresponding to the eigenvector $\lambda_i$. This can be done since the eigenvectors are all linearly independent \cite{Roman2008}. Once the projectors are written, then,
	$$(D-z)^{-1} = \frac{1}{\lambda_1-z} \cdot P_1 + \frac{1}{\lambda_2-z} \cdot P_2 + \frac{1}{\lambda_3-z} \cdot P_3.$$

To find the projection matrices, we must find matrices that have the following properties:

\begin{enumerate}
	\item map each eigenspace to itself, $P_i(E_i) = E_i$;
	\item $P_i^2 = P_i$;
	\item $P_i \cdot v = 0$ when $v \in \underset{j \neq i}{\oplus} Range(P_j)$; and
	\item $P_1 + P_2 + P_3 = I$.
\end{enumerate}

First we consider $P_1$ and $P_3$, both of which are $(m-1)^2 \times (m-1)^2$ matrices:
$$P_1 = \frac{1}{m(m-1)} \left(\begin{array}{cccc}
	1 & \cdots & 1 & m\\
	\vdots & & \vdots & \vdots\\
	1 & \cdots & 1 & m
\end{array}\right)_{(m-1)^2 \times (m-1)^2},$$
$$P_3 = \frac{1}{m(m-2)} \left(\begin{array}{cccc}
	J & \cdots & J & 0\\
	\vdots & & \vdots & \vdots\\
	J & \cdots & J & 0\\
	0 & \cdots & 0 & 0
\end{array}\right)_{(m-1)^2 \times (m-1)^2}$$
where $J$ is the $m \times m$ matrix
$$J = \left(\begin{array}{cccc}
	m-1 & -1 & \cdots & -1\\
	-1 & m-1 & \cdots & -1\\
	\vdots & & \ddots & \vdots\\
	-1 & -1 & \cdots & m-1
\end{array}\right)_{m \times m}.$$

Then we set $P_2 = I - (P_1 + P_3)$ to satisfy the fourth requirement for the projection matrices listed above. While adding $P_1$ and $P_3$, we obtain a new $m \times m$ matrix, $K$. Subtracting $P_1 + P_3$ from $I$, we have:
$$P2 = \left(\begin{array}{ccccc}
	I-K & -K & \cdots & -K & \frac{-1}{m-1}\\
	-K & I-K & \cdots & -K & \frac{-1}{m-1}\\
	\vdots & & \ddots & & \vdots\\
	-K & -K & \cdots & I-K & \frac{-1}{m-1}\\
	\frac{-1}{m(m-1)} & \frac{-1}{m(m-1)} & \cdots & \frac{-1}{m(m-1)} & \frac{m-2}{m-1}
\end{array}\right)_{(m-1)^2 \times (m-1)^2}$$ 
where
$$K = \text{\small{$\frac{1}{m(m-1)(m-2)}$}} \left(\begin{array}{cccc}
	m^2-m-1 & -1 & \cdots & -1\\
	-1 & m^2-m-1 & \cdots & -1\\
	\vdots & & \ddots & \vdots\\
	-1 & -1 & \cdots & m^2-m-1
\end{array}\right)_{m \times m}.$$
It is not difficult to show that these are projecting onto the correct subspaces by checking their action on the eigenvectors.

Now that we have our projection matrices and their corresponding eigenvalues, and since matrix multiplication is distributive, we can write the Schur Complement as: $S(z) = (A-z) - \sum\limits_{i=1}^{3} \frac{1}{\lambda_i-z}BP_iC$. The $1,1$ and $1,2$ entries are the only ones we need to calculate $\phi(z)$ and $R(z)$. The relevant entries from the summands are:
$$\begin{array}{c|ccc}
	& P_1 & P_2 & P3\\
	\hline
	\frac{1}{\lambda_i-z}BP_iC_{1,1} & \frac{1}{m(1-z(m-1))} & 0 & \frac{m-2}{m(2-z(m-1))}\\
	\frac{1}{\lambda_i-z}BP_iC_{1,2} & \frac{1}{m(1-z(m-1))} & 0 & \frac{-(m-2)}{m(m-1)(2-z(m-1))}
\end{array}.$$

Since $(B(D-z)^{-1}C)_{1,1} = \sum\limits_{i=1}^{3} \frac{1}{\lambda_i-z}(BP_iC)_{1,1}$ and $(B(D-z)^{-1}C)_{1,2} = \sum\limits_{i=1}^{3} \frac{1}{\lambda_i-z}(BP_iC)_{1,2}$, we have
	$$S_{1,1} = (1-z) - \frac{m-z(m-1)^2}{m(1-z(m-1))(2-z(m-1))}$$ and
	$$S_{1,2} = -\frac{m-z(m-1)}{m(m-1)(1-z(m-1))(2-z(m-1))}.$$
\end{proof}

\begin{corollary}\label{c:phi}
For any $m \geqslant 3$, the function $\phi(z) = -(m-1)S_{1,2}$ is:
	$$\phi(z) = \frac{m-(m-1)z}{m(2-(m-1)z)(1-(m-1)z)}.$$
\end{corollary}

\begin{corollary}\label{c:R}
For any $m \geqslant 3$, the function $R(z) = 1-\frac{S_{1,1}}{\phi(z)}$ is:
	$$R(z) = 2mz - m(m-1)z^2.$$
\end{corollary}

\begin{figure}[htbp]
	\centering
		\includegraphics[height=100pt]{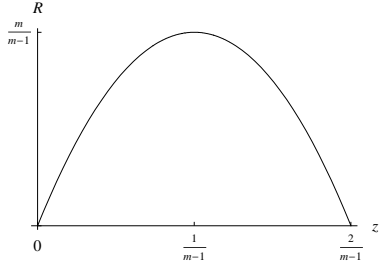}
	\caption{The graph of $R(z)$ for the $m$-branch tree}
	\label{fig:Rm}
\end{figure}

\begin{theorem}\label{thm:full}
For any $m \geqslant 3$, $n \geqslant 0$, and $1 \leqslant k \leqslant n$, the depth-$n$ Laplacian matrix will have eigenvalues such that:
$$\begin{array}{l}
	\textup{mult}_n(0) = 1,\\
	\textup{mult}_n(\frac{m}{m-1}) = 1 + (m-2) \cdot m^n,\\
	\textup{mult}_n(R^{-k}(\frac{m}{m-1})) = 1 + (m-2) \cdot m^{n-k}.
\end{array}$$
In the limit as $n \to \infty$, $\sigma(\Delta_m) = \bigcup_{i\ge 0}\left\{R^{-i}(\frac{m}{m-1})\right\} \cup \{0\} \subseteq [0,\frac{2}{m-1}) \cup \{\frac{m}{m-1}\}$, where the inverse images of $\frac{m}{m-1}$ accumulate to the Julia set of $R(z)$.
\end{theorem}
\begin{proof}
We first define the exceptional set, $E(M_{m,0},M_m) = \sigma(D) \cup \{z:\phi(z) = 0\}$, as in Section \ref{sec:spec-dec}. Since $\phi(z) = 0$ only when $z = \frac{m}{m-1}$, we have
	$$E(M_{m,0},M_m) = \left\{\frac{m}{m-1},\frac{2}{m-1},\frac{1}{m-1}\right\}.$$
We can find the multiplicities of these exceptional values using \Pro.

For $z = \frac{m}{m-1}$, which is in $\sigma(D)$, is not a pole of $\phi(z)$, and $\phi(z) = 0$, we use \pro5 to find its multiplicities:
$$\begin{array}{l}
	\textup{mult}_n(\frac{m}{m-1}) = 1 + (m-2) \cdot m^n.
\end{array}$$

For $z = \frac{2}{m-1}$, which is in $\sigma(D)$, is a pole of $\phi(z)$, and $R'(z) \neq 0$, we use \pro3 to find its multiplicities:
$$\begin{array}{l}
	\textup{mult}_n(\frac{2}{m-1}) = 0.\\
\end{array}$$

For $z = \frac{1}{m-1}$, which is in $\sigma(D)$, is a pole of $\phi(z)$, and $R'(z) = 0$, we use \pro6 to find its multiplicities. Note that $\textup{mult}_0(\frac{1}{m-1}) = 0$ so:
$$\begin{array}{l}
	\textup{mult}_n(\frac{1}{m-1})  = \textup{mult}_{n-1}(\frac{m}{m-1})= 1 + (m-2) \cdot m^{n-1}.\\
\end{array}$$

When $z = 0$, since $z \notin E(M_{m,0},M_m)$, we use \pro1 to find its multiplicities. Given $\textup{mult}_0(0) = 1$ and $R(0) = 0$,
	$$\textup{mult}_n(0) = 1.$$

Since $R^{-1}(0) = \{0, \frac{2}{m-1}\}$ and $\textup{mult}_n(\frac{2}{m-1}) = 0$, there are no additional ancestors from 0. However, as $n$ increases, one must consider $R^{-k}(\frac{m}{m-1})$ for $1 \leqslant k \leqslant n$. Since $R^{-k}(\frac{m}{m-1}) \notin E(M_0,M)$ when $k > 1$, we use \pro1 to determine that
	$$\textup{mult}_n \left(R^{-k} \left(\tfrac{m}{m-1}\right)\right)
	= \textup{mult}_{n-k} \left(\tfrac{m}{m-1}\right)
	= 1 + (m-2) \cdot m^{n-k}.$$

As $n \to \infty$, the only ancestor from $z=0$ is itself. However, since $R(z)$ has zeros at $z = \{0,\frac{2}{m-1}\}$ and maximum at $R(\frac{1}{m-1}) = \frac{m}{m-1}$, $R^{-1}(z) \in [0,\frac{2}{m-1}]$ so long as $z \in [0,\frac{m}{m-1}]$, as seen in Figure \ref{fig:Rm}. Since $\bigcup_{i\ge 0}\left\{R^{-i}\left(\frac{m}{m-1}\right)\right\} \subseteq (0,\frac{2}{m-1}) \cup \{\frac{m}{m-1}\}$, we have
	$$\sigma(\Delta_m) = \bigcup_{i\ge 0}\left\{R^{-i}(\tfrac{m}{m-1})\right\} \cup \{0\} \subseteq [0,\tfrac{2}{m-1}) \cup \{\tfrac{m}{m-1}\}.$$
	
	As before, the inverse images accumulate to the Julia set of $R(z)$ as discussed in \cite{Falconer1990}.
\end{proof}

Table \ref{tab:ancestorm} shows the ancestor-offspring structure of the eigenvalues of the $m$-branch tree for the first few Laplacian matrices. In the table, we denote $R^{-1}(\frac{1}{m-1})=\psi$ and $R^{-1}(\psi) = *$. Blank entries denote exceptional values and if they have non-zero multiplicity they are added at the right hand end of the table.

	\global\newcount\NN
	\def\newz#1{\global\NN=0 $z\in\sigma(M_{#1})$\vbox to 4ex{}}
	\def\newm#1{\global\NN=1 $\textup{mult}_{#1}(z)$\vbox to 4ex{}}
	\def\m#1#2#3{\multicolumn{#1}{c|}{\text{\hskip-.2em\ifcase\NN
	\mathversion{bold}$#2$\or\bfseries\itshape#3\fi\hskip-.2em}}}

\begin{table}[htbp]
	\def\newz#1{\global\NN=0 \text{\small{\hskip-.3em}$z{\in}\sigma(M_{#1})$\vbox to 4ex{}{\hskip-.5em}}}
	\def\newm#1{\small\global\NN=1 {\hskip-.3em}$\textup{mult}_{#1}(z)${\hskip-.5em}\vbox to 4ex{}}
	\def\m#1#2#3{\multicolumn{#1}{c|}{\text{\small\hskip-.5em\ifcase\NN
	\small$#2$\or\bfseries\itshape\small#3\fi\hskip-.5em}}}%
        \centering\small
\begin{tabular}{|l||*{22}{c}}
	\cline{1-9}
		\newz0 &\m401 &\m4{\frac{m}{m-1}}{$m-1$}\\
		\newm0 &\m401 &\m4{\frac{m}{m-1}}{$m-1$}\\
		[2ex]\cline{1-13}
		\newz1 &\m301 &\m1{\text{\tiny{$\frac{2}{m-1}$}}}{} &\m4{\frac{1}{m-1}}{$m-1$} &\m4{\frac{m}{m-1}}{$(m-1)^2$}\\
		\newm1 &\m301 &\m1{\text{\tiny{$\frac{2}{m-1}$}}}{} &\m4{\frac{1}{m-1}}{$m-1$} &\m4{\frac{m}{m-1}}{$(m-1)^2$}\\
		[2ex]\cline{1-4}\cline{6-15}
		\newz2 &\m201 &\m1{\text{\tiny{$\frac{2}{m-1}$}}}{} &\m1{}{} &\m2{\psi_1}{$m-1$} &\m2{\psi_2}{$m-1$}
			&\m4{\frac{1}{m-1}}{$(m-1)^2$} &\m2{\frac{m}{m-1}}{\text{\tiny{$m^3-2m^2+1$}}}\\
		\newm2 &\m201 &\m1{\text{\tiny{$\frac{2}{m-1}$}}}{} &\m1{}{} &\m2{\psi_1}{$m-1$} &\m2{\psi_2}{$m-1$}
			&\m4{\frac{1}{m-1}}{$(m-1)^2$} &\m2{\frac{m}{m-1}}{\text{\tiny{$m^3-2m^2+1$}}}\\
		[2ex]\cline{1-3}\cline{6-16}
		\newz3 &\m101 &\m1{\text{\tiny{$\frac{2}{m-1}$}}}{} &\m1{}{} &\m1{}{}
			&\m1{*}{\text{\tiny{$m-1$}}} &\m1{*}{\text{\tiny{$m-1$}}}	&\m1{*}{\text{\tiny{$m-1$}}} &\m1{*}{\text{\tiny{$m-1$}}}
			&\m2{\psi_1}{\text{\tiny{$(m-1)^2$}}}	&\m2{\psi_2}{\text{\tiny{$(m-1)^2$}}}
			&\m2{\frac{1}{m-1}}{\text{\tiny{$m^3-2m^2+1$}}}	&\m1{\frac{m}{m-1}}{\text{\tiny{$m^4-2m^3+1$}}}\\
		\newm3 &\m101 &\m1{\text{\tiny{$\frac{2}{m-1}$}}}{} &\m1{}{} &\m1{}{}
			&\m1{*}{\text{\tiny{$m-1$}}} &\m1{*}{\text{\tiny{$m-1$}}}	&\m1{*}{\text{\tiny{$m-1$}}} &\m1{*}{\text{\tiny{$m-1$}}}
			&\m2{\psi_1}{\text{\tiny{$(m-1)^2$}}}	&\m2{\psi_2}{\text{\tiny{$(m-1)^2$}}}
			&\m2{\frac{1}{m-1}}{\text{\tiny{$m^3-2m^2+1$}}}	&\m1{\frac{m}{m-1}}{\text{\tiny{$m^4-2m^3+1$}}}\\
		[2ex]\cline{1-16}
\end{tabular}\vskip2ex
	\caption{Ancestor-offspring structure of the eigenvalues on the $m$-branch tree}
	\label{tab:ancestorm}
\end{table}

\begin{proposition}\label{prop:full}
The eigenvalues and multiplicities established in Theorem \ref{thm:full} represent all the eigenvalues of $M_{m,n}$. That is, at any level $n$, the number of eigenvalues is equal to $\textup{dim}_n$.
\end{proposition}
\begin{proof}
We can conclude from Theorem \ref{thm:full} and Table \ref{tab:ancestorm} that all eigenvalues for any level $n$ in $F_m$ lie in either $[0,\frac{2}{m-1})$ or $\{\frac{m}{m-1}\}$. The number of eigenvalues in the depth-$n$ Laplacian of the $m$-branch tree will be
\begin{eqnarray*}
	& & \textup{mult}_n(0) + \textup{mult}_n(\tfrac{m}{m-1}) + \sum\limits_{i=1}^n \textup{mult}_n \left(R^{-i}(\tfrac{m}{m-1})\right)\\
	&=& 1 + (1+(m-2)m^n) + \sum\limits_{i=1}^n (2^{i-1}(1+(m-2)m^{n-i}))\\
	&=& 1 + (1+(m-2)m^n) + (2^n-1) + (m^n-2^n)\\
	&=& 1 + (m-1)m^n = \textup{dim}_n.\\
\end{eqnarray*}
The equation holds because $R^{-1}(z)$ will yield exactly two values for all $z \in (0,\frac{2}{m-1})$. Since $\textup{dim}_n = 1 + (m-1)m^n$, we have found every eigenvalue for the depth-$n$ Laplacian matrix.
\end{proof}

It is of some interest as well to have some idea of the geometric properties of the fractal as well as its vibrational behavior. To this end we need more detailed information about the generating iterated function system than the qualitative description given in the original description of the fractals. 

\begin{proposition}\label{prop:contraction}
For any $n \geqslant 0$ and $m \geqslant 3$, the common contraction factor $c_m$ for all functions in the IFS as measured in the effective resistance metric is:
	$$c_m = \frac{1}{2}.$$
\end{proposition}

\begin{proof} The idea behind this is that if two elements of $V_0$ are separated by resistance one, then those same points viewed in $V_1$ still should have resistance one between them. So you pare off the uninvolved branches and get that they are now at opposite ends of a chain of two copies of $V_0$ each of which would have equal resistance by symmetry so they must have resistance one half. So the distance between neighboring points goes down by a factor of one half each time the level of approximation is increased. For a more formal treatment of these types of calculations see \cite{Kigami1995,Strichartz2006}.
\end{proof}

Once the contraction factor is known it is a direct corollary to calculate the dimension of a given $m-$branch tree. These fractals satisfy the open set condition which is a prerequisite of the formula in this proposition.

\begin{proposition}\label{c:h-dim}
Since the $m$-Branch Trees satisfy the open set condition and are generated by $m$ contraction mappings the Housdorff dimension $s_m$ is given by:
	$$s_m = \frac{-2\log(m)}{\log\left(\frac{1}{2} \right)}.$$
\end{proposition}

\begin{proof} 
A formula for the Hausdorff dimension of a fractal generated by an IFS with a finite number of functions is given in \cite{Falconer1990}. The formula is:$$\sum_{i=1}^m c_i^{s_m} = 1.$$
In our case $c_i = \frac12$ so the calculations are simple and yield $$s_m = \frac{-2\log(m)}{\log(\frac12)}.$$ 
\end{proof}

 This confirms the reasonable intuition that as the number of branches grows so does the dimension of the space. 

\section{Towards the Infinite-Branch Tree}\label{sec:infinity}
With the full spectrum of the $m$-branch tree established in Theorem \ref{thm:full}, we can observe what happens as the number of branches goes towards infinity. From Propositions \ref{p-mult} and \ref{prop:dim} and Theorem \ref{thm:full}, it is clear that the multiplicities of eigenvalues and $\textup{dim}_{m,n}$ both go to infinity. However, we can observe what occurs to the density of the spectra as $m \to \infty$.

As a tool to describe what happens to the spectra of the $\Delta_m$ as $m$ grows we introduce a measure on the complex plane supported on $\sigma(\Delta_m)$. For a $z \in R^{-k}\left(\frac{m}{m-1}\right)$ then for some $n$ large enough $\textup{mult}_n(z)$ is positive. Also $0$ is an eigenvalue with $\textup{mult}_n(0) = 1$ for all $n$. Points, $z$, of these two forms are dense in $\sigma(\Delta_m)$ so we define the measure as follows 
$$\kappa_m = \lim_{n \rightarrow \infty} \frac{1}{dim_n} \sum_{z \in \sigma(\Delta_m)} \textup{mult}_n(z) \delta_z.$$
Now we consider the weak limit of these measures and the support of the limit measure.

\begin{theorem}\label{thm:infinity}
As  $\lim_{m \rightarrow \infty} \sigma(\Delta_m) = \{0,1\}$ in the Hausdorff metric. Also $\lim_{m \rightarrow \infty} \kappa_m = \delta_1$ where the limit is in the weak topology. 
\end{theorem}

\begin{proof}
From Theorem \ref{thm:full}, we have that $$\bigcup_{i\ge 0}\left\{R^{-i}\left(\frac{m}{m-1}\right)\right\} \subseteq \left[0,\frac{2}{m-1}\right] \cup \left\{\frac{m}{m-1}\right\}$$
and 
$$\sigma(\Delta_m) =\bigcup_{i\ge 0}\left\{R^{-i}(\tfrac{m}{m-1})\right\} \cup \{0\} \subseteq \left[0,\tfrac{2}{m-1}\right] \cup \left\{\tfrac{m}{m-1}\right\}.$$
It is easy to see that the only points in the complex plane about any open ball will intersect infinitely many of the $\sigma(\Delta_m)$ are the points $\{0\}$ and $\{1\}$. 

Since $\textup{mult}_n(\frac{m}{m-1}) = 1 + (m-2)m^n$ and $\textup{dim}_n = 1 + (m-1)m^n$, the density $\mathcal{F}$ of $\frac{m}{m-1}$ in the spectrum of $M_n$ as $m \to \infty$ is $$\mathcal{F}\left(\frac{m}{m-1}\right) = \frac{1+(m-2)m^n}{1+(m-1)m^n} = \frac{m^{n+1}-2m^n+1}{m^{n+1}-m^n+1} \rightarrow_{m \rightarrow \infty} 1.$$ And similarly the density of all other eigenvalues which are the ones tending towards zero is $O(\frac{1}{m})$ as $m \rightarrow \infty$.

\end{proof}


\end{document}